\documentclass{amsart}
\usepackage{amssymb, verbatim}
\usepackage{graphicx}

\newcommand{\Lap}{\Delta}

\newcommand{\abs}[1]{{\left\lvert{#1}\right\rvert}}
\newcommand{\norm}[1]{{\left\lVert{#1}\right\rVert}}

\newcommand{\ang}[1]{{\left\langle{#1}\right\rangle}}
\newcommand{\pa}{{\partial}}
\newcommand{\nd}{{\partial_N}}

\newcommand{\ep}{{\epsilon}}

\theoremstyle{plain}
\newtheorem{theorem}{Theorem}[section]

\newtheorem{lemma}[theorem]{Lemma}

\theoremstyle{definition}

\theoremstyle{remark}

\renewcommand\div{\operatorname{div}}

\title{Spreading of quasimodes in the Bunimovich stadium}
\author{Nicolas Burq}
\address{Universit\'e Paris Sud, Math\'ematiques, B\^at 425, 91405 Orsay
  Cedex FRANCE and Institut Universitaire de France}
\email{nicolas.burq@math.u-psud.fr}
\author{Andrew Hassell}
\address{Department of Mathematics, Australian National University, Canberra 0200 ACT, \phantom{ar} AUSTRALIA}
\email{hassell@maths.anu.edu.au}
\author{Jared Wunsch}
\address{Department of Mathematics, Northwestern University,
  Evanston IL 60208, USA}
\email{jwunsch@math.northwestern.edu}
\thanks{This research was partially supported by a Discovery Grant from the
  Australian Research Council (AH), National Science Foundation grants
  DMS-0323021 and DMS-0401323 (JW), and the Institut Universitaire de
  France (NB).  NB and JW gratefully acknowledge the hospitality of the
  Mathematical Sciences Institute of the Australian National University.}

\begin{document}

\begin{abstract} We consider Dirichlet eigenfunctions
$u_\lambda$ of the Bunimovich stadium $S$, satisfying $(\Delta - \lambda^2)
  u_\lambda = 0$.  Write $S = R \cup W$ where $R$ is the central rectangle
  and $W$ denotes the ``wings,'' i.e. the two semicircular regions. It is a
  topic of current interest in quantum theory to know whether
  eigenfunctions can concentrate in $R$ as $\lambda \to \infty$.  We obtain
  a lower bound $C \lambda^{-2}$ on the $L^2$ mass of $u_\lambda$ in $W$,
  assuming that $u_\lambda$ itself is $L^2$-normalized; in other words, the
  $L^2$ norm of $u_\lambda$ is controlled by $\lambda^2$ times the $L^2$
  norm in $W$. Moreover, if $u_\lambda$ is a $o(\lambda^{-2})$ quasimode,
  the same result holds, while for a $o(1)$ quasimode we prove that $L^2$
  norm of $u_\lambda$ is controlled by $\lambda^4$ times the $L^2$ norm in
  $W$.  We also show that the $L^2$ norm of $u_\lambda$ may be controlled
  by the integral of $w \abs{\pa_N u}^2$ along $\partial S \cap W$, where
  $w$ is a smooth factor on $W$ vanishing at $R \cap W$.  These results
  complement recent work of Burq-Zworski which shows that the $L^2$ norm of
  $u_\lambda$ is controlled by the $L^2$ norm in any pair of strips
  contained in $R$, but adjacent to $W$.
\end{abstract}

\maketitle

\section{Introduction}
The Bunimovich stadium $S$ is a planar domain given by the union of a
rectangle $R = \{ (x,y) \mid x \in [-\alpha,\alpha], \ y \in [-\beta,\beta]
\}$ with two ``wings,'' i.e. the two semicircular regions centered at
$(\pm\alpha,0)$ with radius $\beta$ which lie outside $R$.  Geodesic flow
in $S$ (obeying the law of reflection at the boundary) was proved to be
ergodic by Bunimovich \cite{Bun}. It follows from this and from
results of Schnirelman \cite{S}, Zelditch \cite{Z} and Colin de Verdi\`ere
\cite{CdV} that the stadium is quantum ergodic. This means that there is a
density one sequence of Dirichlet eigenfunctions which becomes uniformly
distributed; in particular, along this density one sequence the weak limit
of the $L^2$ mass distribution becomes uniform. One can ask whether the
entire sequence of eigenfunctions becomes uniformly distributed; if so, the
domain is called quantum unique ergodic (QUE). It has been conjectured by
Rudnick and Sarnak that complete surfaces with negative curvature which are
classically ergodic are QUE; this has been proved recently by Lindenstrauss
\cite{L} for arithmetic surfaces.\footnote{With one slight caveat, that the
eigenfunctions are also eigenfunctions of the Hecke operators.}  The
Bunimovich stadium, by contrast, is generally believed to be non-QUE; it is
thought that there is a sequence of eigenfunctions that concentrates in the
rectangle $R$. Little is currently understood about the way in which such
eigenfunctions would concentrate, however.  For example their
(hypothetical) rate of decay outside $R$ is unclear.  The result in the
present paper is intended to shed some light on this question: we show that
any sequence of eigenfunctions (or quasimodes) cannot concentrate very
rapidly inside $R$, by obtaining lower bounds (tending to zero as $\lambda
\to \infty$, but only polynomially) on the $L^2$ mass inside the wings
$W_\pm$.

Let $\Lap =
-\pa_x^2-\pa_y^2$ denote the (nonnegative) Laplacian on $S$ with Dirichlet boundary conditions.  We denote by 
$\norm{\cdot}$ the norm in $L^2(S)$, and by 
$\nd g$ the
outward pointing normal derivative of $g$ at $\pa S.$ We consider a $o(1)$ Dirichlet quasimode  $u_\lambda$ for $\Delta$, by which we mean that we have a sequence 
$\lambda = \lambda_k\to \infty$ of real numbers and a corresponding sequence
$u_{\lambda} \in H^2(S)$ satisfying
\begin{equation}\label{quasimode}
\begin{aligned}
(\Lap-\lambda^2) u_{\lambda} &= f_{\lambda},\\
 u_{\lambda}\mid_{\partial S} &= 0,\\
\norm{u_\lambda} &=1,
\end{aligned}
\end{equation}
where
\begin{equation}\label{qm-order}
\norm{f_\lambda} = o(1) \text{ as }\lambda \to \infty.
\end{equation}
We more generally define a $O(\lambda^{-j})$ or $o(\lambda^{-j})$
quasimode by modifying the right-hand side of \eqref{qm-order}
accordingly. Of course a sequence of eigenfunctions is a $o(\lambda^{-j})$ quasimode for any $j$. 

It is easy to see that a $O(1)$ quasimode can be localized to a small
rectangle of the form $[\gamma, \delta] \times [-\beta, \beta]$, where
$[\gamma, \delta]$ is an arbitrary subinterval of $[-\alpha, \alpha]$;
indeed the family $u_{(n+1/2)\pi/\beta} = \phi(x) \cos((n+1/2) \pi y
/\beta)$ is (after normalization) such a quasimode, where $\phi$ is any
nonzero smooth function supported in $[\gamma, \delta]$. By constrast, an
$o(1)$ quasimode cannot be so localized: Burq-Zworski \cite{BZ} have shown
that the $L^2$ norm of $u_\lambda$ is controlled by (that is, bounded above
by a constant times) its $L^2$ norm in the union of any two rectangles of
the form $([-\alpha, \gamma_1] \times [-\beta, \beta]) \cup ([\gamma_2,
\alpha] \times [-\beta, \beta])$. In particular, for a $o(1)$ quasimode,
the $L^2$ mass cannot shrink to a closed region disjoint from the wings of
the stadium as $\lambda \to \infty$.

Although the stadium is classically ergodic, there is a codimension one
invariant set for the classical flow, consisting of vertical ``bouncing
ball'' orbits parallel to the $y$-axis and within the rectangle $R$, and
the union of these orbits is the most likely place where localization of
eigenfunctions, or more generally $o(1)$ quasimodes, can occur\footnote{The
explicit quasimode in the paragraph above concentrates along a subset of
these orbits}. There is a rather convincing plausibility argument in the
physics literature due to Heller and O'Connor \cite{HOC} which indicates
that a density-zero sequence of eigenfunctions, with eigenvalues
$((n+1/2)\pi/\beta)^2 + O(1)$, does concentrate to some extent at these
bouncing-ball orbits. The rigorous essence of this argument has been
developed by Donnelly \cite{D} who showed that there are sequences of
functions lying in the range of spectral projectors $E_{I_n}(\Delta)$,
where $I_n$ are intervals of the form $[((n+1/2)\pi/\beta)^2 -C,
((n+1/2)\pi/\beta)^2 + C]$ which concentrate at the bouncing-ball
orbits.\footnote{This was shown for surfaces without boundary containing a
flat cylinder, but the arguments go through for the stadium.} On the other
hand, the result of Burq-Zworski \cite{BZ} shows that such localization
cannot be too extreme: the control region must extend to the boundary of
the rectangle.

Our main result here is that we may in fact push our control region outside
the rectangle altogether and into the wings, in return for a loss, either
from restriction to the boundary, or of powers of $\lambda.$ To state this concisely, it is
convenient to introduce an auxiliary coordinate in the wings given by
$w=\abs{x}-\alpha;$ thus $w$ is nonnegative on the wings and vanishes
exactly on the vertical lines $R \cap W.$

\begin{theorem}\label{ourtheorem}
There is a $C > 0$, depending only on $\alpha/\beta$, 
such that 
any family $u_\lambda$ satisfying \eqref{quasimode} obeys the estimates 
\begin{equation}
\label{normderiv}
\|f_\lambda\|^2 +\int_{\pa S \cap W} \!\! w \, \abs{\nd u_\lambda}^2 \, dl \geq C \ ,
\end{equation}
\begin{equation}
\label{L2}
\|f_\lambda\|^2+ \lambda^8 \norm{u_\lambda}_{L^2(W)}^2 \geq C
\end{equation}
and
\begin{equation}
\label{L2bis}
\lambda^2 \|f_\lambda\|^2 + \lambda^4 \norm{u_\lambda}_{L^2(W)}^2 \geq C.
\end{equation}

Therefore, if $u_\lambda$ is a $o(1)$ quasimode, we have for sufficiently large $\lambda$
\begin{equation}\begin{aligned}
\| w^{\frac1{2}} \partial_N u \|_{L^2(\pa S \cap W)} & \geq C,\\ 
\norm{u_\lambda}_{L^2(W)} &\geq C \lambda^{-4},
\label{o1}
\end{aligned}
\end{equation}
while if $u_\lambda$ is a $o(\lambda^{-2})$ quasimode (e.g.\ an eigenfunction),
\begin{equation}\label{betterquasimode}
\norm{u_\lambda}_{L^2(W)} \geq C \lambda^{-2}
\end{equation}
\end{theorem}

Note that the results of Theorem~\ref{ourtheorem} still leave open the
possibility of quasimodes concentrated along bouncing-ball orbits in the
rectangle with $o(1)$ mass in the wings.  They also do not rule out the
possibility that all the energy in the wings may asymptotically concentrate
in a boundary layer near $R.$


\section{Preliminaries to $L^2$ estimates}
Our main tool is positive commutator estimates, which we use in the following form:
\begin{lemma} \label{lemma:rellich}
Let $u$ be real, equal to zero at $\pa S,$ and satisfy $(\Delta -
\lambda^2) u = f,$ where $f$ is smooth. Then for any real vector field $A$,
\begin{equation}
 \langle u, [\Delta - \lambda^2,A]u \rangle  = 
 \langle (2 Au+ (\div A) u, f \rangle +
\int_{\pa S} (\partial_N u) Au \, dl. 
\label{Rellich}\end{equation}
\end{lemma}

\begin{proof} 
We integrate twice by parts, using the Dirichlet boundary conditions in the
first instance, to write 
$$ \ang{u, [\Lap-\lambda^2, A] u} = \ang{f,A u} + \int_{\pa S} \pa_N u Au
\, dl - \ang{u, A f}.
$$
Applying Green's Theorem to the last term now gives two terms: $\ang{Au, f}+ \ang{(\div
  A) u, f}.$  Since $u$ and $f$ are real this yields the desired identity.
\end{proof}

We also record here an inequality that will be of use in
estimating derivative terms.
\begin{lemma}\label{lemma:gradientestimate}
Let $u,$ $f$ be as in Lemma~\ref{lemma:rellich}.

Then for all $s>0,$ for $\lambda$ sufficiently large,
$$
\norm{\nabla u}^2 \leq  C_s (\lambda^{\max(2,s)}\norm{u}^2 + \lambda^{-s} \norm{f}^2).
$$
\end{lemma}

\begin{proof}
We compute
\begin{align*}
\norm{\nabla u}^2 &= \int_S u_x^2+u_y^2 \, dA\\ &= \int_S (\Lap u ) u \,
dA\\ &= \lambda^2 \int_S u^2 \,dA + \ang{f,u}.
\end{align*}
Applying Cauchy-Schwarz to  $\ang{f,u}$ gives the estimate.
\end{proof}


\section{Proof of \eqref{normderiv}}

It suffices to prove \eqref{normderiv} under the assumption that $u_\lambda$, and hence $f_\lambda$, are real, since we can treat the real and imaginary parts separately. We make this assumption from now on. 

We begin with the standard commutator $[\Lap,x\pa_x] = -2 \pa_x^2.$ 
Applying \eqref{Rellich} with $A = x \pa_x$, and dropping the subscript on $u_\lambda$, we have
\begin{equation}\begin{aligned}
\ang{u_x, u_x} = -\ang{\pa_x^2 u,u}&= \ang{[\Lap-\lambda^2, x \pa_x] u,u} 
\\
&= \int_{\pa S} x\pa_x u \, \nd u\, dl
+ \int (2 x \pa_xu + u) f \, dA ;
\end{aligned}\label{boundary}\end{equation}
in the last equation we integrated twice by parts, using for a second time
the fact that $(\Lap-\lambda^2)u=f$ as well as the fact that $u$ satisfies
Dirichlet boundary conditions, hence integration by parts produces
boundary terms only where derivatives land on both factors of $u.$

Now at every boundary point we may decompose $x \pa x$ into $p \pa_l+q \nd$
where $\pa_l$ is differentiation tangent to the boundary.  Of course
$\pa_l$ annihilates $u.$  Now since $\pa_x$ is tangent to the upper and
bottom sides of the rectangle we find that the boundary integral in 
\eqref{boundary} is only over $\pa S \cap W.$  Moreover, as $\pa_x$ is
tangent to the top and bottom of the circles
forming the boundaries of the wings, we have $q = O(w)$ on $\pa S \cap
W.$  Hence we have shown that
\begin{equation}
\| u_x \|^2 \leq \int_{\pa S \cap W} O(w)\abs{\nd u}^2 \, dl + \ep
\int (u^2 + u_x^2) \, dA + C \int f^2 \, dA.
\label{useful}\end{equation}
We may  absorb the $\epsilon \| u_x \|^2$ term, then apply the Poincar\'e inequality, and absorb the $\epsilon \| u \|^2$ to obtain
\begin{equation}
\norm{u}_{L^2(S)}^2 \leq C \int_{\pa S \cap W} w \abs{\nd u}^2 \, dl+ C  \norm{f}^2
\label{ndd}\end{equation}
which is the first part of our theorem, as we took $u$ to be
$L^2$-normalized.

\section{Proof of \eqref{L2}}
To prove this estimate we start from 
\begin{equation}\label{normderivbis}
\| u_x \|^2 \leq C \int_{\partial S} w_+ | \nd u|^2 \, dl +C \norm{f}^2
\end{equation}
which follows directly from the considerations of the previous section, and estimate the boundary integral term. We 
shall obtain upper
bounds of the form
$$
\lambda^8 \int_W u^2\, dA + \| f \|^2,
$$
and
$$
\lambda^4 \int_W u^2\, dA + \lambda^2 \| f \|^2,
$$
thus proving \eqref{L2} and \eqref{L2bis}.

We shall perform this estimate in three separate regions in the wing.  
Region I is the near-rectangular region, in a boundary layer where
$w \leq \delta\lambda^{-2}.$  Region II will be outside the boundary layer, 
where  $\delta\lambda^{-2} \leq w \leq \beta/2.$  Region III will be the far outer region $w \geq \beta/2.$

\begin{figure}
\includegraphics[scale=.3]{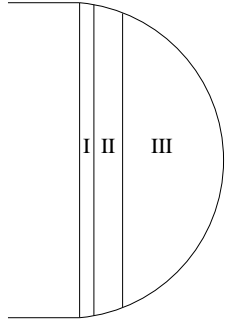}
\caption{The three regions of interest in $W.$}
\end{figure}

We begin with Region III, far away from the rectangle. In this case we
employ Lemma~\ref{lemma:rellich} where $A$ is the operator $\phi(x)
\partial_x$, where $\phi$ is supported where $w > \beta/4$, say,  
equal to $1$ where $w > \beta/2$, and with $\partial_x \phi \geq 0$. Then \eqref{Rellich} gives us, with $P =
[\Delta, A]$,
$$
\int_{\partial S} \phi \partial_x u \, \partial_N u \leq \abs{\langle Pu, u
\rangle}+\abs{\ang{\phi_x f,u}} +2\ang{\phi u_x,f}.
$$
Note that $P = -2 \phi_x \partial^2_{xx} - \phi_{xx}
\partial_x$. Thus  the LHS is bounded ($\forall \epsilon>0$) by
$$
\abs{\ang{- \phi_x u_{xx}, u}} + \ep \norm{u_x}^2 + C (\norm{u}_{L^2(W)}^2+\norm{f}^2)
$$ 
(where of course $C$ depends on $\ep$).
We can add the positive term $\int_S \phi_x (\partial_y u)^2$ to this
estimate. Integrating by parts in $y$ gives us
$$
\abs{\ang{ \phi_x (-u_{xx}- u_{yy}), u}} + \ep \norm{u_x}^2 + C (\norm{u}_{L^2(W)}^2+\norm{f}^2).
$$  Using the positivity of the integrand, we thus obtain an estimate
\begin{equation}
\abs{\int_{\partial S\cap \text{III}} \phi  \, \partial_x u \, \partial_N u\, dl} \leq C \lambda^2
\norm{u}_{L^2(W)}^2 + C \norm{f}^2+ \ep \norm{u_x}^2
\label{bdy-est-1}\end{equation}
with $C$ depending on $\ep>0.$

Now we work on Region I, within a $O(\lambda^{-2})$ boundary layer along
the rectangle.  We again apply Lemma~\ref{lemma:rellich}, this time with
$A= x\pa_x+y\pa_y.$  Since $A$ is a tangential vector plus a positive
multiple of $\pa_N$ all along $\pa S$ we obtain
$$ \int_{\partial S} (\partial_N u)^2\, dl \leq \abs{ \langle u, [\Delta-
  \lambda^2, A]u\rangle } + 2 \abs{ \langle u, f \rangle }+ 2 \abs{\ang{Au,f}}.
$$ 
Using $[\Delta - \lambda^2, A] = 2\Delta$ and Cauchy-Schwarz this becomes 
$$
\int_{\partial S} \abs{\partial_N u}^2 \, dl \leq C \lambda^2   \| u \|^2
+ \norm{Au}^2 + C\| f \|^2.
$$  Restricting to $\pa S \cap
\text{I}$ in the integrand, we can estimate $w$ in $L^\infty$ by $\delta
\lambda^{-2}$, and this gives 
\begin{equation} \int_{\pa S
\cap \text{I}} w_+ \abs{\partial_N u}^2 \, dl \leq \delta C (\|u\|^2+ \lambda^{-2}
(\norm{f}^2+\norm{Au}^2)).
\end{equation}
Using Lemma~\ref{lemma:gradientestimate}, we may estimate $\norm{Au}^2$ by
$C (\lambda^2\norm{u}^2+\norm{f}^2)$. Hence we may finally write
\begin{equation}
\label{bdy-est-3}
 \int_{\partial S\cap \text{I}} w_+ \abs{\partial_N u}^2 \, dl \leq \delta C_0 \| u \|^2 +
\delta C \lambda^{-2}\norm{f}^2;
\end{equation}
note that in the above construction, $C_0$ can in fact be chosen independent of $\delta.$

Finally, we estimate in Region II.  To begin with, we note that for $w \geq \delta \lambda^{-2}$, we can estimate $w_+$ by $ \delta^{-1} \lambda^2
w_+^2$,  so we have
$$
\int_{\pa S \cap \text{II}} w_+ \abs{\pa_N u}^2 \, dl
\leq C\delta^{-1}\int_{\pa S \cap \text{II}} \lambda^2 w_+^2 \,  \chi(y)\pa_y u  \, \pa_N u \, dl;
$$
 here we take $\chi$ supported in $|y| > \beta/20$, and equal to $-1$ for
$y < -\beta/10$ and $+1$ for $y > \beta/10$, so that $\chi(y)\pa_y$ is a
positive multiple of $\pa_N$ plus a tangential component on $\pa S \cap
\text{II}.$ To estimate further, we employ Lemma~\ref{lemma:rellich} with
the commutant $$A = \lambda^2 w_+^2 \chi(y) \partial_y.$$ The point of this
commutant is that we have given ourselves two powers of $w$ which will
``absorb'' two integrations by parts in $x$ \emph{without any boundary terms
at $w=0$}; this is crucial since we know of no way to deal with such boundary terms. (On the other hand, we pay the price of additional powers of $\lambda$ with this gambit.)
Thus we obtain, setting $Q = [\Delta, A]$, 
\begin{equation}\label{region2}
\int_{\partial S} \lambda^2 w_+^2 \chi(y) \partial_y u \,  \partial_N u\, dA \leq \Big| 
\langle Qu, u \rangle \Big| +  \lambda^2 \abs{ \int_S w_+^2 \chi_y(y) u f}
+ 2 \abs{\ang{Au, f}}
\end{equation}
We can estimate the second term on the RHS by $\lambda^4 \| u \|_{L^2(W)}^2
+ C\| f \|^2$.

Now consider the terms involving the operator $Q$. This is given by
$$
Q = \lambda^2 \Big( -4 w_+ \chi \partial_x \partial_y -2 w_+^2 \chi_y \partial^2_{yy} - 2 H(w) \chi  \partial_y - w_+^2\chi_{yy} \partial_y \Big)  
$$
where $H(\cdot)$ is the Heaviside step-function, $H(w) = 0$ for $w < 0$ and $1$ for $w \geq 0$. 
To treat the terms involving one derivative, e.g. the third term above, we integrate by parts:
$$
-2\lambda^2 \int_S H(w) \chi(y) \,  \partial_y u \,  u =  \lambda^2 \int_S H(w) \chi_y(y) \, u^2 
$$
which is therefore estimated by $\lambda^2 \| u \|_{L^2(W)}^2$. The fourth term is estimated in exactly the same way. 

Thus we are left to estimate
 \begin{equation}
 \lambda^2  \Big( -4 \langle w_+ \chi(y) \partial^2_{xy} u,  u \rangle
-2  \langle w_+^2 \chi_y(y) \partial^2_{yy} u,  u \rangle \Big)
\end{equation}
Integrating the first term by parts in $x$ and the second term by parts in $y$ gives us two principal terms 
 \begin{equation}\label{foobar}
4 \lambda^2  \Big( \langle w_+^{1/2} \partial_x u, w_+^{1/2} \partial_y u \rangle
+ \langle w_+^2 \chi_y \partial_y u, \partial_y u \rangle \Big)
\end{equation}
together with two other terms
 \begin{equation*} 4 \lambda^2  \Big( \langle H(w) \partial_y u,  u \rangle
+ \langle w_+^2 \chi_{yy}(y) \partial_y u,  u \rangle \Big)
\end{equation*}
which are estimated in the same way as the first order terms above. 
We apply Cauchy-Schwarz to the first term in \eqref{foobar}, while in the second term, which is positive, we replace $w_+^2$ by $w_+$ (which is larger, up to a
constant multiple) and then integrate by
parts again, getting (up to another first-order term estimated as above) an upper bound for \eqref{foobar} of the form
\begin{equation}\label{region2gradient}
C \lambda^2 \int_S w_+ \abs{\partial_x u}^2 +  w_+ \abs{\partial_y u}^2
\, dA.
\end{equation}
Now we integrate by parts again, getting
\begin{equation}
C \lambda^2 \int_S \Big( w_+ (-u_{xx} - u_{yy}) u + H(w) u \, u_x \Big)\, dA.
\label{eee}\end{equation}
Writing $-u_{xx} - u_{yy} = \lambda^2 u + f$, we can estimate the integrand of \eqref{eee} by 
\begin{multline*}
\lambda^4 w_+ \abs{u}^2 + \lambda^2 w_+ f u + \lambda^2 H(w) u u_x 
\\ \leq \lambda^4 w_+ \abs{u}^2 + \frac1{2} \big( \lambda^4 w_+  u^2 + w_+ f^2 \big) + C \lambda^4 H(w) u^2 + \epsilon u_x^2.
\end{multline*}
Hence we may estimate \eqref{foobar}, by as
$$
C \lambda^4
\norm{u}_{L^2(W)}^2 + \ep \norm{u_x}^2+ C \norm{f}^2 .
$$

The term $2\abs{\ang{Au,f}}$
is bounded in a similar manner:
By Cauchy-Schwarz we may estimate it by
$$
C \norm{f}^2 + C\lambda^4\norm{w_+ u_y}^2
$$
or by
$$
C \lambda^2 \norm{f}^2 + C\lambda^{2}\norm{w_+ u_y}^2,
$$
as we prefer.  Our estimate for \eqref{region2gradient} turns the latter
estimate into
$$
C \lambda^2 \norm{f}^2 +C \lambda^4
\norm{u}_{L^2(W)}^2 + \ep \norm{u_x}^2.
$$
On the other hand, treating $\lambda^4\norm{w_+ u_y}^2$ in the same manner
gives a bound by
$$
C \norm{f}^2 +C \lambda^8
\norm{u}_{L^2(W)}^2 + \ep \norm{u_x}^2.
$$

On the support of $\chi$ and at the boundary, $\chi(y) u_y$ is a
positive multiple of $\partial_N u$. So the upshot is that \eqref{region2} now gives
\begin{equation}
\int_{\partial S} \lambda^2 w_+^2 \chi(y) (\partial_N u)^2 \leq 
C \lambda^8
\norm{u}_{L^2(W)}^2 + \ep \norm{u_x}^2 + C \norm{f}^2.
\label{bdy-est-2}\end{equation}
and
\begin{equation}
\int_{\partial S} \lambda^2 w_+^2 \chi(y) (\partial_N u)^2 \leq 
C \lambda^4
\norm{u}_{L^2(W)}^2 + \ep \norm{u_x}^2 + C \lambda^2 \norm{f}^2.
\label{bdy-est-2bis}\end{equation}

At last we can estimate the boundary integral term in \eqref{normderivbis} by a combination of
\eqref{bdy-est-1}, \eqref{bdy-est-3}, and
\eqref{bdy-est-2}/\eqref{bdy-est-2bis}, obtaining
$$
\int_{\partial S} w_+ |\nd u|^2 \, dl \leq \delta C_0\|u\|^2+ C \lambda^8
\norm{u}_{L^2(W)}^2 + 2 \ep \norm{u_x}^2 + C \norm{f}^2.
$$
and
$$
\int_{\partial S} w_+ |\nd u |^2 \, dl \leq \delta C_0\|u\|^2+ C \lambda^4
\norm{u}_{L^2(W)}^2 + 2 \ep \norm{u_x}^2 + C \lambda^2 \norm{f}^2.
$$
 Here $C$ depends on $\ep,$ $\delta$ but $C_0$ is independent of
$\delta.$ We now combine this with \eqref{normderivbis}. 
Absorbing the
$\norm{u_x}^2$ and $\norm{u}^2$ terms on the LHS by taking $\delta$ and $\ep$
sufficiently small, we obtain \eqref{L2}, \eqref{L2bis}.

\section{Concluding remarks}

The estimates presented here are certainly not optimal. The powers of
$\lambda$ appearing in Theorem~\ref{ourtheorem} can probably be improved
using refinements of the methods used here, but it seems unlikely that one
could achieve an optimal result with them. We have not, therefore, attempted to
obtain the best possible powers of $\lambda$, but have rather tried to
present a poynomial lower bound on $\| u_\lambda \|_{L^2(W)}$ with a simple
proof.

It would be of great interest to obtain a polynomial lower bound on the
$L^2$ mass of $u_\lambda$ in a subregion of $W$ which is a positive
distance from $R$ (i.e. region III of the previous section). We do not know
whether such a bound holds, but it does not seem to be obtainable using the
methods of this paper; possibly it might yield to the use of more
sophisticated tools from microlocal analysis.

\end{document}